\documentclass[12pt]{amsart}

\usepackage{color}
\usepackage{amsmath, amsthm, amssymb}
\usepackage{amsfonts}
\usepackage[ansinew]{inputenc}
\usepackage[dvips]{epsfig}
\usepackage{graphicx}
\usepackage[english]{babel}
\usepackage{hyperref}
\theoremstyle{plain}
\newtheorem{thm}{Theorem}[section]
\newtheorem{cor}[thm]{Corollary}

\theoremstyle{definition}

\theoremstyle{remark}

\numberwithin{equation}{section}

\begin{document}

\title{Fractional Laplacian: Pohozaev identity and nonexistence results}

\author{Xavier Ros-Oton}

\address{Universitat Polit\`ecnica de Catalunya, Departament de Matem\`{a}tica  Aplicada I, Diagonal 647, 08028 Barcelona, Spain}
\email{xavier.ros.oton@upc.edu}

\thanks{The authors were supported by grants MTM2011-27739-C04-01 (Spain) and 2009SGR345 (Catalunya)
}

\author{Joaquim Serra}

\address{Universitat Polit\`ecnica de Catalunya, Departament de Matem\`{a}tica  Aplicada I, Diagonal 647, 08028 Barcelona, Spain}

\email{joaquim.serra@upc.edu}

\maketitle

\begin{abstract}

\vskip 0.5\baselineskip

In this note we present the Pohozaev identity for the fractional Laplacian.
As a consequence of this identity, we prove the nonexistence of nontrivial bounded solutions to semilinear problems with supercritical nonlinearities in star-shaped domains.

\vskip 1\baselineskip

\noindent{\sc R\'esum\'e.} 
Dans cette note, nous pr\' esentons l'identit\'e de Pohozaev pour le Laplacien fractionnaire.
Comme cons\'equence de cette identit\'e, nous prouvons la non-existence de solutions non triviales pour les probl\`emes semi-lin\'eaires avec nonlin\'earit\'e sur-critique
dans des domaines \'etoil\'es.
\end{abstract}

\vspace{6mm}

\section{Introduction}
\label{intro}

Let $s\in(0,1)$, and consider the fractional elliptic problem
\begin{equation}\label{eq}\left\{ \begin{array}{rcll} (-\Delta)^s u &=&f(u)&\textrm{ in }\Omega \\
u&=&0&\textrm{ in }\mathbb R^n\backslash\Omega,\end{array}\right.\end{equation}
in a bounded domain $\Omega\subset\mathbb R^n$, where
\[(-\Delta)^s u(x)=c_{n,s}\mbox{PV}\int_{\mathbb R^n}\frac{u(x)-u(y)}{|x-y|^{n+2s}}dy\]
and $c_{n,s}$ is a normalization constant.

When $s=1$, a celebrated result of S. I. Pohozaev states that any solution of (\ref{eq}) satisfies an identity, which is known as the Pohozaev identity \cite{P}. This classical result yields, as an immediate consequence, the nonexistence of nontrivial bounded solutions to (\ref{eq}) for supercritical nonlinearities $f$ in star-shaped domains $\Omega$.
In this note we present the fractional version of this identity, that is, a generalization of the Pohozaev identity which applies to problem (\ref{eq}).
This result will be proved in a forthcoming paper \cite{RS}, and reads as follows.
Here, since the solution $u$ is bounded, the notions of energy and viscosity solution agree (see \cite{RS}).


\begin{thm}\label{pohz} Let $\Omega$ be a bounded and $C^{1,1}$ domain, $f$ be a locally Lipschitz function, $u\in H^s(\mathbb R^n)$ be a bounded solution of (\ref{eq}), and
$\delta(x)={\rm dist}(x,\partial\Omega)$.
Then
\[u/\delta^s|_{\Omega}\in C^{0,\alpha}(\overline\Omega)\qquad \textrm{for some }\ \alpha\in(0,1),\]
meaning that $u/\delta^s|_{\Omega}$ has a continuous extension to $\overline\Omega$ which is
$C^{0,\alpha}(\overline{\Omega})$, and the following identity holds
\[(2s-n)\int_{\Omega}uf(u)dx+2n\int_\Omega F(u)dx=\Gamma(1+s)^2\int_{\partial\Omega}\left(\frac{u}{\delta^s}\right)^2(x\cdot\nu)d\sigma,\]
where $F(t)=\int_0^tf$, $\nu$ is the unit outward normal to $\partial\Omega$ at $x$, and $\Gamma$ is the Gamma function.
\end{thm}


Note that the function $u/\delta^s|_{\partial\Omega}$ plays the role that $\partial u/\partial\nu$ plays in the classical Pohozaev identity.
Moreover, our Pohozaev identity for $s=1$ is obviously the classical one, since $u/\delta|_{\partial\Omega}=\partial u/\partial\nu$ and $\Gamma(2)=1$.


As an immediate consequence of this identity we improve some recent results of M.M. Fall and T. Weth \cite{FW} on nonexistence of solutions to problem (\ref{eq}) with supercritical nonlinearities in star-shaped domains.


\begin{cor}\label{cornonexistence} Let $\Omega$ be a bounded, $C^{1,1}$, and star-shaped domain, and $f$ be a locally Lipschitz function.~If
\begin{equation}\label{supercritic} \frac{n-2s}{2n}uf(u)\geq\int_0^u f(t)dt\qquad\textrm{for all}\ \ u\in\mathbb R, \end{equation}
then problem (\ref{eq}) admits no positive bounded solution.

Moreover, if the inequality in (\ref{supercritic}) is strict, then (\ref{eq}) admits no nontrivial bounded solution.
\end{cor}


The proof of the nonexistence results in \cite{FW} uses the method of moving spheres and, therefore, in \cite{FW} solutions are assumed to be positive.
Our nonexistence result is the first one allowing changing-sign solutions.
As in \cite{FW}, we may allow also nonlinearities $f(x,u)$ depending on $x\in\Omega$; see \cite{RS}.


In addition to Theorem \ref{pohz}, in \cite{RS} we will also obtain the following integration by parts formula.


\begin{thm}\label{intparts} Let $\Omega$ be a bounded and $C^{1,1}$ domain. Assume that $u$ and $v$ are bounded solutions of $(-\Delta)^s u=g(x,u)$ and $(-\Delta)^sv=h(x,v)$ in $\Omega$, and $u\equiv v\equiv0$ in $\mathbb R^n\backslash\Omega$, for some $g,h\in C^{0,1}_{\rm loc}(\overline\Omega\times\mathbb R)$, and let $\delta(x)={\rm dist}(x,\partial\Omega)$.
Then, $u/\delta^s|_\Omega$ and $v/\delta^s|_\Omega$ have $C^{0,\alpha}(\overline\Omega)$ extensions and it holds
\[\int_\Omega (-\Delta)^su\ v_{x_i}dx=-\int_\Omega u_{x_i}(-\Delta)^sv\ dx+\Gamma(1+s)^2\int_{\partial\Omega}\left(\frac{u}{\delta^{s}}\right)\left(\frac{v}{\delta^{s}}\right) \nu_i\ d\sigma\]
for all $i\in\{1,...,n\}$, where $\nu$ is the unit outward normal to $\partial\Omega$ at $x$ and $\Gamma$ is the Gamma function.
\end{thm}

\section{Sketch of the proofs}

Let us next give an sketch of the proof of our fractional Pohozaev identity in star-shaped domains.
The identity in non star-shaped domains is deduced afterwards using a partition of the unity, as shown in \cite{RS}.
The main idea of the proof is to use
\[u_\lambda(x)=u(\lambda x), \quad \lambda>1,\]
as a test function in the weak formulation of problem (\ref{eq}) and then differentiate the obtained identity with respect to $\lambda$ at $\lambda=1$.
However, this apparently simple formal procedure requires a quite involved analysis when it is put into practice.
Namely, it needs fine boundary regularity results for $u$ and $u/\delta^s$, (the main one obtained through a Krylov boundary Harnack method), as well as the precise behavior of $(-\Delta)^{s/2}u$ in all of $\mathbb R^n$. Recall that we denote $\delta(x)={\rm dist}(x,\partial\Omega)$.

Although Corollary \ref{cornonexistence} follows immediately from Theorem \ref{pohz}, we give here a short proof of the nonexistence result for supercritical nonlinearities $f$, that is,
when inequality in (\ref{supercritic}) is strict. The proof of this result follows the same method that we use to establish the Pohozaev identity, but it does not require the precise analysis mentioned above.

The proofs start by showing that
\begin{equation}\label{first}\int_\Omega(x\cdot \nabla u) (-\Delta)^su\ dx=\left.\frac{d}{d\lambda}\right|_{\lambda=1^+}\int_\Omega u_\lambda (-\Delta)^su\ dx.\end{equation}
This equality follows from the gradient estimate
$|\nabla u|\leq C\delta^{s-1}$ in $\Omega$ and the dominated convergence theorem.
This gradient estimate is proved using standard regularity arguments, detailed in \cite{RS}.

Now we use $u_\lambda$, $\lambda>1$, as a test function for problem (\ref{eq}).
At this point it is crucial to assume that the domain $\Omega$ is star-shaped,
which guarantees that $u_\lambda\equiv0$ in $\mathbb R^n\backslash \Omega$.
We obtain
\[\int_\Omega u_\lambda(-\Delta)^{s}u\ dx = \int_{\mathbb R^n}(-\Delta)^{s/2}u_{\lambda}(-\Delta)^{s/2}u\ dx=
\lambda^s\int_{\mathbb R^n}w_\lambda w\, dx.\]
where $w(x)=(-\Delta)^{s/2}u(x)$ and $w_{\lambda}(x)=w(\lambda x)$.

Moreover, with the change of variables $x\mapsto \lambda^{-1/2}x$ this integral becomes
\[\int_\Omega u_\lambda(-\Delta)^su\ dx=\lambda^s\int_{\mathbb R^n}w_{\lambda}w\ dx=\lambda^{\frac{2s-n}{2}}\int_{\mathbb R^n}w_{\sqrt{\lambda}}w_{1/\sqrt{\lambda}}\,dx,\]
which leads to
\begin{eqnarray}
\int_\Omega(\nabla u\cdot x)(-\Delta)^sudx&=&\left.\frac{d}{d\lambda}\right|_{\lambda=1^+}\left\{\lambda^{\frac{2s-n}{2}}\int_{\mathbb R^n}w_{\sqrt{\lambda}}w_{1/\sqrt{\lambda}}dx\right\}\nonumber \\
&&\hspace{-20mm}=\frac{2s-n}{2}\int_{\mathbb R^n}|(-\Delta)^{s/2}u|^2dx +\left.\frac{d}{d\lambda}\right|_{\lambda=1^+}\int_{\mathbb R^n}w_{\sqrt{\lambda}}w_{1/\sqrt{\lambda}}\,dx\nonumber \\
&&\hspace{-20mm}=\frac{2s-n}{2}\int_{\Omega}u(-\Delta)^s u\ dx+ \frac 12\left.\frac{d}{d\lambda}\right|_{\lambda=1^+}\int_{\mathbb R^n}w_{{\lambda}}w_{1/{\lambda}}\,dx\label{cadira}
.\end{eqnarray}

Furthermore, since $(-\Delta)^su=f(u)$ in $\Omega$ and 
\[\int_\Omega(x\cdot\nabla u)(-\Delta)^su\, dx=\int_\Omega x\cdot\nabla F(u)dx=-n\int_\Omega F(u)dx,\]
(\ref{cadira}) reads as
\begin{equation}\label{casipoh}
-n\int_\Omega F(u)dx=\frac{2s-n}{2}\int_\Omega uf(u)dx+\frac 12\left.\frac{d}{d\lambda}\right|_{\lambda=1^+}\int_{\mathbb R^n}w_{{\lambda}}w_{1/{\lambda}}dx.
\end{equation}
Thus, the Pohozaev identity is equivalent to
\begin{equation}\label{derivadaIlambda}\left.\frac{d}{d\lambda}\right|_{\lambda=1^+}
I_\lambda=-\Gamma(1+s)^2\int_{\partial\Omega}\left(\frac{u}{\delta^s}\right)^2(x\cdot\nu)d\sigma,
\end{equation}
where
\[ I_\lambda=\int_{\mathbb R^n}w_{{\lambda}}w_{1/{\lambda}}dx.\]
This equality is the difficult part of the proof of the Pohozaev identity.

The quantity $\left.\frac{d}{d\lambda}\right|_{\lambda=1^+}\int_{\mathbb R^n}w_{\lambda}w_{1/\lambda}$ vanishes for any $C^1(\mathbb R^n)$ function $w$ ---as can be seen by differentiating under the integral sign.
Instead, the function $w=(-\Delta)^{s/2} u$ has a singularity along $\partial\Omega$, and a crucial part of our proof consists of establishing
the precise behavior of $(-\Delta)^{s/2}u$ near $\partial\Omega$ (from both inside and outside $\Omega$), namely \[(-\Delta)^{s/2}u(x)=c_1\left\{\log\delta(x)+c_2\chi_{\Omega}(x)\right\}\frac{u}{\delta^s}(x^*)+h(x),\]
where $c_1$ and $c_2$ are constants, $x^*$ is the nearest point to $x$ on $\partial\Omega$, and $h$ is a $C^{0,\alpha}$ function.

In contrast with equality (\ref{derivadaIlambda}) ---which is delicate to prove---, the inequality
\begin{equation}\label{ineq} \left.\frac{d}{d\lambda}\right|_{\lambda=1^+}I_\lambda\leq0\end{equation}
follows easily from Cauchy-Schwarz.
Indeed, we have
\[I_\lambda\leq \|w_\lambda\|_{L^2}\|w_{1/\lambda}\|_{L^2}=\|w\|_{L^2}^2=I_1,\]
and thus $\frac{I_\lambda-I_1}{\lambda-1}\leq0$, which yields (\ref{ineq}).
Finally, by (\ref{casipoh}) and (\ref{ineq}) we deduce
\[n\int_\Omega F(u)dx\geq \frac{n-2s}{2}\int_{\Omega}uf(u)dx.\]
This yields the nonexistence result in star-shaped domains for supercritical nonlinearities, that is, Corollary \ref{cornonexistence} with strict inequality in (\ref{supercritic}).


\section*{Acknowledgements}

The authors thank Xavier Cabr\'e for his help an all his comments on this note.


\end{document}